\documentclass[11pt]{article}
\usepackage{}
\usepackage[total={6in, 9in}]{geometry}
 \usepackage{amsfonts}
\usepackage{amsthm}
\usepackage{amssymb}
\usepackage{amsmath,amsfonts}
\usepackage{mathrsfs}
\usepackage{cases}
\usepackage{latexsym,bm}
\usepackage{indentfirst}
\usepackage{color}
\usepackage{ifpdf}
\usepackage{graphicx}
\usepackage{psfrag}
\usepackage{enumerate}
\usepackage{dsfont}
\usepackage[
pdfauthor={shan},
pdftitle={Ore-type condition  for hamiltonicity in tough graphs},
pdfstartview=XYZ,
bookmarks=true,
colorlinks=true,
linkcolor=blue,
urlcolor=blue,
citecolor=blue,
bookmarks=true,
linktocpage=true,
hyperindex=true
]{hyperref}

\newtheorem{THM}{\textbf{Theorem}}
\newtheorem{THMs}{\textbf{Theorem}}[section]

\newtheorem{LEM}{\textbf{Lemma}}
\newtheorem{CLA}{\textbf{Claim}}
\newtheorem{CON}{\textbf{Conjecture}}

\newcommand{\pf}{\noindent\textbf{Proof}.\quad}

\newcommand{\iC}{\overset{\leftharpoonup }{C}}

\newcommand{\oC}{\overset{\rightharpoonup }{C}}

\newcommand{\CC}{\mathcal{C}}

\def\dist{{\fam0 dist}}

\begin{document}
\title{An Ore-type condition  for hamiltonicity in tough graphs}
\author{ Songling Shan \\ 
	\medskip  Illinois State  University, Normal, IL 61790\\
	\medskip 
	{\tt sshan12@ilstu.edu}
}

\date{\today}
\maketitle

\emph{\textbf{Abstract}.}
Let $G$ be a $t$-tough graph on $n\ge 3$ vertices for some $t>0$.  It was  shown by Bauer et al. in 1995 that if the minimum degree of $G$ is greater than $\frac{n}{t+1}-1$, then $G$ is hamiltonian. 
In terms of Ore-type 
hamiltonicity conditions, the problem was only studied when $t$ is between 1 and 2. 
 In this paper, we show
that if the degree sum of any two nonadjacent vertices of $G$ is greater than $\frac{2n}{t+1}+t-2$, then $G$ is hamiltonian.


\vspace{2mm}

\section{Introduction}

\section{Introduction}

We consider only simple  graphs. 
Let $G$ be a graph.
Denote by $V(G)$ and  $E(G)$ the vertex set and edge set of $G$,
respectively. Let $v\in V(G)$, $S\subseteq V(G)$, and $H\subseteq G$. 
Then  $N_G(v)$   denotes the set of neighbors
of $v$ in $G$, $d_G(v):=|N_G(v)|$ is the degree of $v$ in $G$, 
and $\delta(G):=\min\{d_G(v): v\in V(G)\}$ is the minimum degree of  $G$. 
Define 
$\deg_G(v,H)=|N_G(v)\cap V(H)|$, $N_G(S)=(\bigcup_{x\in S}N_G(x))\setminus S$, 
and  we write $N_G(H)$ for $N_G(V(H))$. 
Let $N_H(v)=N_G(v)\cap V(H)$
and $N_H(S)=N_G(S)\cap V(H)$. 
We use $G[S]$ and $G-S$ to denote the subgraphs of $G$ induced by  $S$ and $V(G)\setminus S$, respectively. 
For notational simplicity we write $G-x$ for $G-\{x\}$.
Let $V_1,
V_2\subseteq V(G)$ be two disjoint vertex sets. Then $E_G(V_1,V_2)$ is the set
of edges in $G$  with one end in $V_1$ and the other end in $V_2$. For two integers $a$ and $b$, let $[a,b]=\{i\in \mathbb{Z}\,:\,   a\le i \le b\}$.

Throughout this paper,  if not specified, 
we will assume $t$ to be a nonnegative real number. The number of components of  a graph $G$ is denoted by $c(G)$. 
The graph $G$ is said to be {\it $t$-tough\/} if $|S|\ge t\cdot
c(G-S)$ for each $S\subseteq V(G)$ with $c(G-S)\ge 2$. The {\it
	toughness $\tau(G)$\/} is the largest real number $t$ for which $G$ is
$t$-tough, or is  $\infty$ if $G$ is complete. This concept, a
measure of graph connectivity and ``resilience'' under removal of
vertices, was introduced by Chv\'atal~\cite{chvatal-tough-c} in 1973.
It is  easy to see that  if $G$ has a hamiltonian cycle
then $G$ is 1-tough. Conversely,
Chv\'atal~\cite{chvatal-tough-c}
conjectured that
there exists a constant $t_0$ such that every
$t_0$-tough graph is hamiltonian.
Bauer, Broersma and Veldman~\cite{Tough-counterE} have constructed
$t$-tough graphs that are not hamiltonian for all $t < \frac{9}{4}$, so
$t_0$ must be at least $\frac{9}{4}$ if Chv\'atal's toughness conjecture is true.

Chv\'atal's toughness conjecture  has
been verified when restricted to a number of graph
classes~\cite{Bauer2006},
including planar graphs, claw-free graphs, co-comparability graphs, and
chordal graphs.  In general, the conjecture is still wide open. 
In finding hamiltonian cycles in graphs, sufficient  conditions such as Dirac-type  and 
Ore-type conditions are the most classical ones. 

\begin{THMs}[Dirac's Theorem~\cite{Dirac-theorem}]
	If $G$ is a graph on $n\ge 3$  vertices with $\delta(G) \ge \frac{n}{2}$, then $G$
	is hamiltonian. 
\end{THMs} 

Define $\sigma_2(G)=\min\{d_G(u)+d_G(v)\,:\,  u,v\in V(G),\,\, \text{$u$ and $v$ are nonadjacent}\}$ 
if $G$ is noncomplete, and define $\sigma_2(G)=\infty$ otherwise. Ore's Theorem, as a generalization of 
Dirac's Theorem, is stated below. 

\begin{THMs}[Ore's Theorem~\cite{Ore-Theorem}]\label{thm:ore}
	If $G$ is a graph on $n\ge 3$ vertices with $\sigma_2(G) \ge n$, then $G$
	is hamiltonian. 
\end{THMs}

Analogous to Dirac's Theorem, Bauer, Broersma, van den Heuvel, and Veldman~\cite{MR1336668} proved the following result 
by incorporating the toughness of the graph. 

\begin{THMs}[Bauer et al.~\cite{MR1336668}]\label{degree-tough}
	Let $G$ be a $t$-tough graph on $n\ge 3$ vertices.  If $\delta(G) > \frac{n}{t+1}-1$, then 
	$G$ is hamiltonian. 
\end{THMs}

A natural question here is whether we can  find an Ore-type condition  involving the toughness of $G$
that generalizes Theorem~\ref{degree-tough}.  Various theorems were proved prior to 
Theorem~\ref{degree-tough} by only taking $\tau(G)$ between 1 and 2. 
Jung in 1978~\cite{MR499116}
showed that if $G$ is a 1-tough graph on $n\ge 11$ vertices with $\sigma_2(G) \ge n-4$, then $G$ is hamiltonian.
In 1991, Bauer,  Chen,    and Lasser~\cite{BCL} showed that the degree bound in Jung's Theorem can be slightly lowered if
$\tau(G)>1$. The result states that  
if $G$ is a graph on $n\ge 30$ vertices with  $\tau(G)>1$ and $\sigma_2(G)\ge n-7$, then $G$ is
hamiltonian. In 1989/90, Bauer, Veldman, Morgana, and Schmeichel~\cite{MR1032634} showed that 
if $G$ is a 2-tough graph on $n\ge 3$ vertices with $\sigma_2(G) \ge \frac{2n}{3}$, then $G$ is hamiltonian (a consequence of Corollary 16 from~\cite{MR1032634}).
In this paper, we obtain the following result, which provides an Ore-type condition  involving  $\tau(G)$
that guarantees a hamiltonian cycle in a graph.

\begin{THM}\label{main}
	Let $G$ be a $t$-tough graph on $n\ge 3$ vertices.  If $\sigma_2(G) > \frac{2n}{t+1}+t-2$, then 
	$G$ is hamiltonian. 
\end{THM}
In fact, we believe that the following stronger statement might be true. 
\begin{CON}\label{con}
	Let $G$ be a $t$-tough graph on $n\ge 3$ vertices.  If $\sigma_2(G) > \frac{2n}{t+1}-2$, then 
	$G$ is hamiltonian. 
\end{CON}

Considering both toughness and degree sum conditions such as in Theorem~\ref{main} and Conjecture~\ref{con} is an approach to 
investigate Chv\'atal's toughness conjecture while the conjecture remains open. However, 
in light of the conjecture, those results might only be relevant for some small values of $t$.

For odd integers $n\ge 3$,  the complete bipartite graph $G:=K_{\frac{n-1}{2}, \frac{n+1}{2}}$ is $\frac{n-1}{n+1}$-tough and 
satisfies  $\sigma_2(G)=n-1 =\frac{2n}{1+\frac{n-1}{n+1}}-2$.  However, $G$ is not hamiltonian. 
Thus,  if true, 
the degree sum condition in Conjecture~\ref{con} would be best possible. 
In fact, for odd integers $n\ge 3$, any graph from the family $\mathcal{H}=\{H_{\frac{n-1}{2}} + \overline{K}_{\frac{n+1}{2}}: \text{$H_{\frac{n-1}{2}}$  is any graph on $\frac{n-1}{2}$ vertices}\}$  
is an extremal graph, where ``$+$'' represents the join of two graphs.   
In light of    the results mentioned in the paragraph right above Theorem~\ref{main} and Chv\'atal's toughness conjecture,  it suggests 
that  $t$-tough non-hamiltonian  graphs $G$ with  $\sigma_2(G) = \frac{2n}{t+1}-2$  exist only when  $t<1$. Furthermore, 
by looking at the non-hamiltonian $t$-tough graphs $G$ with $t<1$ and $\delta(G)= \frac{n}{t+1}-1$, which  are exactly the graphs in the family $\mathcal{H}$, 
it suggests that when $t<1$, any non-hamiltonian $t$-tough graph $G$ with $t<1$ and $\sigma_2(G) = \frac{2n}{t+1}-2$ belongs to 
the family $\mathcal{H}$.  So we propose the following conjecture. 
\begin{CON}\label{con2}
	Let $G$ be a $t$-tough graph on $n\ge 3$ vertices.  If $\sigma_2(G) = \frac{2n}{t+1}-2$ and 
	$G$ is non-hamiltonian, then $G\in \mathcal{H}$.  
\end{CON}

In attempting to prove Conjecture~\ref{con} by contradiction,    the most difficult case   to deal with is 
when $G$ has a cycle $C$ of length  $n-1$ and $G-V(C)$ is just a single vertex component $H$.  
It seems very hard to  deduce any nontrivial property of $G$ using the $\sigma_2(G)$ and toughness conditions.  
However, by adding $t$ to the $\sigma_2(G)$  bound,  
vertices in $V(C)\setminus N_C(H)$  can be shown to  have degree bigger than $\frac{n}{t+1}+t-1$. 
This degree condition allows us to find $|N_C(H)|$ disjoint subgraphs each of order $t+2$ 
such that there is no edge between any two of them. Then  we get to use the toughness condition to 
give a smaller upper bound on $|N_C(H)|$ ($|N_C(H)| \le \frac{n}{2(t+1)}-\frac{1}{2}$),  which plays a key role in the proof of Theorem~\ref{main}.
Therefore, it might require a completely different approach to confirm Conjecture~\ref{con}.

The remainder of this paper is organized as follows: in Section 2, we introduce some notation and preliminary 
results,  and in Section 3, we prove Theorem~\ref{main}.

\section{Preliminary results}

Let $G$ be a graph and $\lambda$ a positive integer. Following~\cite{MR696293}, a cycle $C$ of $G$ is 
a  \emph{$D_\lambda$-cycle} if every  component of $G-V(C)$ has order less than $\lambda$. 
Clearly, a $D_1$-cycle is just a hamiltonian cycle.
We denote by $c_\lambda(G)$  the number of components of $G$
with order at least $\lambda$, and write $c_1(G)$ just as $c(G)$. 
Two subgraphs $H_1$ and $H_2$ of $G$ are  \emph{remote} if they are disjoint  and there is no
edge of $G$ joining a vertex of $H_1$ and a vertex of $H_2$. 
For a subgraph $H$ of $G$, let $d_G(H)=|N_G(H)|$ be the degree of  $H$ in $G$. 
We denote by $\delta_\lambda(G)$ the minimum
degree of a connected subgraph of order $\lambda$ in $G$. Again $\delta_1(G)$
is just $\delta(G)$.  

Let $C$ be an oriented cycle, and we assume that the orientation is clockwise throughout the rest of this paper. For $x\in V(C)$,
denote the immediate successor of $x$ on $C$ by $x^+$ and the immediate  predecessor of $x$ on $C$ by $x^-$.
For $u,v\in V(C)$, $u\oC v$  denotes the segment of $C$
starting at $u$, following $C$ in the orientation,  and ending at $v$.
Likewise, $u\iC v$ is the opposite segment of $C$ with endpoints as $u$
and $v$.  Let $\dist_{\oC}(u,v)$  denote the length of the path $u\oC v$.
For any vertex $u\in V(C)$ and any positive integer $k$, define 
$$
L_u^+(k)=\{v\in V(C)\,:\,   \dist_{\oC}(u,v) \in [1,k]\},  \,\, \text{and} \,\, L_u^-(k)=\{v\in V(C)\,:\,   \dist_{\oC}(v,u)  \in [1,k]\}, 
$$
to be the set of $k$ consecutive successors of $u$ and  the set of $k$ consecutive predecessors of $u$, respectively. 
A \emph{chord} of $C$ is an edge $uv$ with $u,v\in V(C)$ and $uv\not\in E(C)$. 
Two chords $ux$ and $vy$ that do not share any endvertices of $C$ are \emph{crossing}  
if the four vertices $u, x, v, y$ appear  along  $\oC$  in the order $u,v, x, y$ or $u, y, x, v$. 
Hereafter, all cycles under consideration are oriented. 

A path $P$ connecting two vertices $u$ and $v$ is called 
a {\it $(u,v)$-path}, and we write $uPv$ or $vPu$ in order to specify the two endvertices of 
$P$. Let $uPv$ and $xQy$ be two paths. If $vx$ is an edge, 
we write $uPvxQy$ as
the concatenation of $P$ and $Q$ through the edge $vx$.

For an integer $\lambda \ge 1$, if a graph $G$ contains a $D_{\lambda+1}$-cycle $C$  but no $D_\lambda$-cycle,
then  $V(G) \setminus V(C) \ne \emptyset$. Furthermore,  $G-V(C)$ has a component  of order $\lambda$. 
The result below with  $d_G(H)$ replaced by $\delta_\lambda(G)$ was proved in~\cite{MR1336668}. 
\begin{LEM}\label{lem:dominating-cycle}
	Let $G$ be a $t$-tough $2$-connected graph of order $n$.
	Suppose $G$ has a $D_{\lambda+1}$-cycle   but no $D_\lambda$-cycle. Let $C$ be a $D_{\lambda+1}$-cycle of $G$ such that $c_\lambda(G-V(C))$ is minimum.  Then $n\ge (t+\lambda) (d_G(H)+1)$
	for any component $H$ of $G-V(C)$ with order $\lambda$.
\end{LEM}

\pf Let $k=d_G(H)$, which equals the total number of neighbors of vertices of $H$ on $C$. 
We assume the $k$ neighbors are $v_1, \ldots, v_k$ and appear in the same order along $\oC$. 
For each $i\in [1,k]$, 
and each $v\in V(v_i^+\oC v_{i+1}^-)$,  where $v_{k+1}:=v_1$, we let $\CC(v)$ be the set of components of $G-V(C)$
that have a vertex joining to $v$ by an edge in $G$. Note that $H\not\in \CC(v)$. 
Let $w_i^*\in V(v_i^+\oC v_{i+1}^-)$ be the vertex with $\dist_{\oC}(v_i,w_i^*)$ minimum  
such that 
$$
\sum\limits_{D\in \bigcup\limits_{v\in V(v_i^+\oC w_i^*)} \CC(v)}|V(D)|+ |V(v_i^+\oC w_i^*)| \ge \lambda. 
$$
If such a vertex $w_i^*$ exists, 
let $L_{v_i}^*(\lambda)$  be the union of the vertex set $V(v_i^+\oC w_i^*)$
and all those vertex sets of graphs in $\bigcup\limits_{v\in V(v_i^+\oC w_i^*)} \CC(v)$;
if such a vertex $w_i^*$ does not exist, let $L_{v_i}^*(\lambda)=L_{v_i}^+(\lambda)$. 
Note that when $w_i^*$ exists, by its definition, $w_i^*\in V(v_i^+\oC v_{i+1}^-)$. Thus  $V(v_i^+\oC w_i^*)\cap V(v_j^+\oC w_j^*) =\emptyset$ if both $w_i^*$ and $w_j^*$ exist for distinct $i,j\in [1,k]$. 

To prove Lemma~\ref{lem:dominating-cycle},  
it suffices to show that $L_{v_1}^*(\lambda), \ldots, L_{v_k}^*(\lambda)$ and $H$ are pairwise remote. 
Since in that case, if we let $S=V(G)\setminus \left( (\bigcup_{i=1}^k L_{v_i}^*(\lambda)) \cup V(H)\right)$,  
then $|S| \le n-(k+1) \lambda$ and $c(G-S) = k+1$. As $G$ is $t$-tough, we get 
$$
n-(k+1) \lambda \ge |S| \ge t\cdot c(G-S) = t(k+1), 
$$
giving $n\ge (t+\lambda)(k+1)$. 

Below, we show that $L_{v_1}^*(\lambda), \ldots, L_{v_k}^*(\lambda)$ and $H$ are pairwise remote. 
It suffices to prove Statement (a):  $\dist_{\oC}(v_i,v_{i+1}) \ge \lambda+1$ if for some $i \in [1,k]$ it holds that $L_{v_i}^*(\lambda)=L_{v_i}^+(\lambda)$, where $v_{k+1}:=v_1$ when $i=k$ (this implies that each $L_{v_i}^*(\lambda)$ and $H$ are remote),  and $L_{v_i}^*(\lambda) \cap L_{v_j}^*(\lambda)=\emptyset$ for every two distinct   $i,j\in [1,k]$; and  Statement (b): 
$E_G(L_{v_i}^*(\lambda), L_{v_j}^*(\lambda))=\emptyset$ for every two distinct   $i,j\in [1,k]$.
Let $v_i^*\in N_H(v_i), v^*_j\in N_H(v_j)$ and $P$ be a $(v_i^*,v^*_j)$-path of $H$.

For Statement (a), 
it suffices to show that if for some $i \in [1,k]$ it holds that $L_{v_i}^*(\lambda)=L_{v_i}^+(\lambda)$, then $\dist_{\oC}(v_i,v_{i+1}) \ge \lambda+1$, where $v_{k+1}:=v_1$ when $i=k$;
and that for distinct $i,j\in [1,k]$, $v\in L_{v_i}^*(\lambda) \cap V(v_i^+ \oC v_{i+1}^-)$ and $u\in L_{v_j}^*(\lambda) \cap V(v_j^+ \oC v_{j+1}^-)$, 
we have $\CC(v)\cap \CC(u) = \emptyset$.  We prove the statement by contradiction. 
If $L_{v_i}^*(\lambda)=L_{v_i}^+(\lambda)$ for some $i\in [1,k]$
but $\dist_{\oC}(v_i,v_{i+1})  \le \lambda$, we then let  
$C^*=v_i\iC v_{i+1}v^*_{i+1}Pv_i^*v_i$. Since $H$ has order $\lambda$ and no vertex of $H$ is adjacent in $G$ to any internal vertex of $v_i\oC v_{i+1}$,   it follows that each component of $H-V(P)$ is a component of $G-V(C^*)$ of order at most $\lambda-1$
and $v_i^+\oC v_{i+1}^-$ is contained in a component of $G-V(C^*)$
with order at most $\lambda-1$ since $L_{v_i}^*(\lambda)=L_{v_i}^+(\lambda)$. Thus  $C^*$ is a $D_{\lambda+1}$-cycle of $G$ with $c_\lambda(G-V(C^*))<c_\lambda(G-V(C))$, contradicting the choice of $C$.  
If for some distinct $i,j\in [1,k]$, $v\in L_{v_i}^*(\lambda) \cap V(v_i^+ \oC v_{i+1}^-)$ and $u\in L_{v_j}^*(\lambda) \cap V(v_j^+ \oC v_{j+1}^-)$, 
we have $\CC(v)\cap \CC(u) \ne \emptyset$,  
we  then further choose $v$ closest  to $v_i$ and $u$ closest to $v_j$  along $\oC$ with the property.  Thus  for any $w_i \in V(v_i^+\oC v^-)$ and any 
$w_j\in V(v_j^+\oC u^-)$,  it holds that $\CC(w_i)\cap \CC(w_j)= \emptyset$.  
Let $D\in \CC(v)\cap \CC(u)$ and  $v', u'\in V(D)$ such that $vv', uu'\in E(G)$. Let $P'$
be a $(v',u')$-path of $D$ and  $C^*=v_iv_i^*Pv_j^*v_j\iC vv'P'u'u\oC v_i$.  
Since $H$ has order $\lambda$ and no vertex of $H$ is adjacent in $G$ to any vertex in $v_i^+\oC v^-$ or any vertex in $v_j^+ \oC u^-$, it follows that each component of $H-V(P)$ is a component of $G-V(C^*)$ of order at most $\lambda-1$. 
Furthermore, by the choices of $v$ and $u$, the components of $G-V(C^*)$
that respectively contain $v_i^+\oC v^-$ and $v_j^+ \oC u^-$ are disjoint. 
Since $V(v_i^+ \oC v^-)$  is a proper subset of $L_{v_i}^*(\lambda) \cap V(v_i^+ \oC v_{i+1}^-)$ and  $V(v_j^+ \oC u^-)$  is a proper subset of $L_{v_j}^*(\lambda) \cap V(v_j^+ \oC v_{j+1}^-)$, it follows by the definitions of $L_{v_i}^*$ and $L_{v_j}^*$
that the components of $G-V(C^*)$
that respectively contain $v_i^+\oC v^-$ and $v_j^+ \oC u^-$  have order at most $\lambda -1$. 
Thus  $C^*$ is a $D_{\lambda+1}$-cycle of $G$ with $c_\lambda(G-V(C^*))<c_\lambda(G-V(C))$, contradicting the choice of $C$.   The argument above verifies Statement (a). 

For Statement (b), 
assume to the contrary that $E_G(L_{v_i}^*(\lambda), L_{v_j}^*(\lambda)) \ne \emptyset$ for some distinct $i,j\in [1,k]$. 
Applying Statement (a),   we know that  $L_{v_i}^*(\lambda) \cap L_{v_j}^*(\lambda)=\emptyset$.
Since there is no edge between any two components of $G-V(C)$,  $E_G(L_{v_i}^*(\lambda), L_{v_j}^*(\lambda)) \ne \emptyset$ implies that there exist $x\in L_{v_i}^*(\lambda) \cap V(v_i^+ \oC v_{i+1}^-)$ and $y\in L_{v_j}^*(\lambda) \cap V(v_j^+ \oC v_{j+1}^-)$ 
such that $xy \in E(G)$. 
We choose  $x\in L_{v_i}^*(\lambda) \cap V(v_i^+ \oC v_{i+1}^-)$ with $\dist_{\oC}(v_i,x)$ minimum
and $y\in L_{v_j}^*(\lambda) \cap V(v_j^+ \oC v_{j+1}^-)$  with $\dist_{\oC}(v_j,y)$ minimum such that 
$xy\in E(G)$. By this choice of $x$ and $y$,  it follows that 
$E_G(V(v_i^+ \oC x^-), V(v_j^+ \oC y^-))=\emptyset$.
Let $C^*=v_i\iC yx\oC v_jv_j^*Pv_i^*v_i$. Since $H$ has order $\lambda$ and 
no vertex of $H$ is adjacent in $G$ to any  vertex of $v_i^+ \oC x^-$ or $v_j^+ \oC y^-$ by the fact that  $v_i^+ \oC x^- \subseteq v_i^+\oC v_{i+1}^-$  and  $v_j^+ \oC y^- \subseteq v_j^+\oC v_{j+1}^-$ from Statement (a),   
it follows that 
each component of   $H-V(P)$ is a component of $G-V(C^*)$ of order at most $\lambda-1$. Also
$v_i^+ \oC x^-$  and $v_j^+ \oC y^-$ are contained in distinct  components of $G-V(C^*)$ 
each of order at most $\lambda -1$.  
Thus  $C^*$ is a $D_{\lambda+1}$-cycle of $G$ with $c_\lambda(G-V(C^*))<c_\lambda(G-V(C))$, contradicting the choice of $C$. 
This verifies Statement (b) and completes the proof of Lemma~\ref{lem:dominating-cycle}. 
\qed

\begin{LEM}[\cite{p2p3}]\label{lem:idependentset-size}
	Let $t>0$ and $G$ be a non-complete $n$-vertex $t$-tough  graph. Then $|W|\le \frac{n}{t+1}$ for every independent set $W$ in $G$. 
\end{LEM}

The following lemma provides a way of extending a cycle $C$  provided that the vertices 
outside $C$ have many neighbors on $C$.  The proof follows from Lemma~\ref{lem:idependentset-size}
and is very similar to the proof of Lemma 10 in~\cite{p2p3}.

\begin{LEM}\label{cycle-extendabilit y}
	Let $t\ge 1$  and $G$ be an $n$-vertex  $t$-tough graph, 
	and let $C$ be a non-hamiltonian cycle of $G$. 
	If $x\in V(G)\setminus V(C)$ satisfies  $\deg_G(x, C)>\frac{n}{t+1}-1$,  then $G$ has a cycle $C'$ such that $V(C')=V(C)\cup \{x\}$. 	
\end{LEM}

\section{Proof of Theorem~\ref{main}}

We may assume that $G$ is not a complete graph. Thus $G$ is $2\lceil t\rceil$-connected as it is $t$-tough. Suppose to the contrary that $G$ is not hamiltonian. 
By Theorem~\ref{degree-tough},  we have  $\delta(G)\le \frac{n}{t+1}-1$. Since $\delta(G) \ge  2\lceil t\rceil$, we get  

\begin{equation}\label{eqn:lower-bound-on-n}
	n\ge (t+1) (2\lceil t\rceil+1). \nonumber
\end{equation}

\begin{CLA}\label{claim:t>=1}
	$t>1$. 
\end{CLA}

\pf Assume to the contrary that $t\le 1$. 
By Ore's result, Theorem~\ref{thm:ore}, and the assumption that $G$ is not hamiltonian, we get  $\sigma_2(G) \le n-1$. 
Thus 
\begin{eqnarray*}
	\frac{2n}{t+1}+t-2< \sigma_2(G)  \le  n-1. 
\end{eqnarray*}
This  gives $\frac{2n}{t+1}+t<n+1$. Let $f(t)=\frac{2n}{t+1}+t$.
Since  $t\le 1$ and $n\ge (t+1) (2\lceil t\rceil+1)>(t+1)^2$, $f'(t)=\frac{(t+1)^2-2n}{(t+1)^2}<0$. Thus the minimum 
value of $f(t)$ is achieved at $t=1$ and $f(1)=n+1$,
showing a contradiction to   $\frac{2n}{t+1}+t<n+1$. 
\qed

Since $t>1$ and $G$ is not complete, $G$
is $2$-connected and so contains cycles. 
We choose $\lambda \ge 0$ to be a smallest integer 
such that $G$ admits no $D_\lambda$-cycle but a $D_{\lambda+1}$-cycle. 
Then we choose   $C$
to be a longest $D_{\lambda+1}$-cycle 
such that $c_\lambda(G-V(C))$ is minimum.  As $G$ is not hamiltonian, 
we have 
$\lambda\ge 1$. 
Thus $V(G)\setminus V(C) \ne \emptyset$. 
Since $\lambda$ is taken to be minimum, $G-V(C)$
has a component  $H$ of order $\lambda$.   
Let 
\begin{equation}\label{eqn:def-of-W}
	W=N_C(V(H)) \quad \text{and} \quad  \omega=|W|.\nonumber     
\end{equation}
Since $G$ is a connected $t$-tough  graph,  it follows that 
$
\omega \ge  2\lceil t\rceil. 
$
On the other hand,  Lemma~\ref{lem:dominating-cycle} implies that 
$
\omega \le \frac{n}{t+ \lambda}-1.
$

\begin{CLA}\label{claim:H-and-W-size}
	$\lambda+ \omega  \le \frac{n}{t+1}$. 
\end{CLA}

\pf  Assume to the contrary that $\lambda+ \omega  > \frac{n}{t+1}$. 
If $\lambda=1$, then  $H$ has only one vertex and $\omega>\frac{n}{t+1}-1$. By Lemma~\ref{cycle-extendabilit y},
we can find a cycle $C'$ with $V(C')=V(C)\cup V(H)$, contradicting the choice of $C$. 
Thus $\lambda \ge 2$. 
Since $ 2t\le \omega \le \frac{n}{t+ \lambda}-1 \le  \frac{n}{t+2}-1$,
we have  $n\ge (t+2)(2t+1)$. 
By Lemma~\ref{lem:dominating-cycle}, we have 
\begin{eqnarray*}
	n &\ge& (\lambda+t)(\omega+1) \\
	&>&  ( \frac{n}{t+1}-\omega+t) (\omega+1)  \quad \text{($\lambda+ \omega  > \frac{n}{t+1}$ by assumption.)}\\
	&\ge & \begin{cases}
		(\frac{n}{t+1}-2t+t)(2t+1), & \text{if $f(\omega)=( \frac{n}{t+1}-\omega+t) (\omega+1)$ is increasing;}\\
		(\frac{n}{t+1}-\frac{n}{t+2}+1+t)\frac{n}{t+2}, & \text{if $f(\omega)=( \frac{n}{t+1}-\omega+t) (\omega+1)$ is decreasing;}
	\end{cases}\\
	& \ge & \begin{cases}
		n+\frac{tn}{t+1}-2t^2-t \ge n+\frac{t(t+2)(2t+1)}{t+1}-2t^2-t>n+t(2t+1)-2t^2-t=n, &\\
		\frac{n}{(t+1)(t+2)} \frac{n}{t+2}+\frac{(t+1)n}{t+2} \ge \frac{n}{(t+1)(t+2)}\frac{(t+2)(2t+1)}{t+2}+\frac{(t+1)n}{t+2}>\frac{n}{t+2}+\frac{(t+1)n}{t+2}=n,
	\end{cases}
\end{eqnarray*}
reaching a contradiction. 
\qed 

\begin{CLA}\label{claim:one-component}
	$H$ is the only component of $G-V(C)$ and $H$ is a complete subgraph of $G$. 
\end{CLA}

\pf Suppose $H^* \ne H$ is another component of $G-V(C)$. 
Since $\sigma_2(G)>\frac{2n}{t+1}+t-2$, Claim~\ref{claim:H-and-W-size}
implies that $|V(H^*)|+|N_C(V(H^*))|>\frac{n}{t+1}+t-1>\frac{n}{t+1}$. 
Repeating exactly the same argument for $|V(H^*)|+|N_C(V(H^*))|$ as in the proof of Claim~\ref{claim:H-and-W-size}
leads to a contradiction. Thus $H$ is the only component of $G-V(C)$. 
Since $\lambda+ \omega  \le \frac{n}{t+1}$
by Claim~\ref{claim:H-and-W-size}  and  $\sigma_2(G)>\frac{2n}{t+1}+t-2$, 
every two distinct   vertices of $H$ are adjacent. 
Thus $H$ is a complete subgraph of $G$. 
\qed 

Since $H$ is the only component of $G-V(C)$,  every vertex $v\in V(C)\setminus W$ is only adjacent in $G$ to vertices on $C$. 
As  $d_G(u) \le \frac{n}{t+1}-1$  for any $u\in V(H)$ by  Claim~\ref{claim:H-and-W-size}, 
using $\sigma_2(G)>\frac{2n}{t+1}+t-2$,   
we have 
\begin{equation}\label{eqn:degree-on-C}
	\deg_G(v,C) >\frac{n}{t+1}+t-1 \quad \text{for any $v\in V(C)\setminus W$}. 
\end{equation}
Equation~\eqref{eqn:degree-on-C} allows us to construct the vertex sets $L_u^+(t+2)$ for each $u\in W$. 
For notation simplicity,  we use $L_u^+$ for $L_u^+(t+2)$.

\begin{CLA}\label{claim:segments}
	For any two distinct vertices $u,v\in W$, $\dist_{\oC}(u,v) \ge t+3$ and 
	$E_G(L_u^+, L_v^+)=\emptyset$. 
\end{CLA}

\pf Let $u^*\in N_H(u), v^*\in N_H(v)$ and $P$ be a $(u^*,v^*)$-path of $H$. 
For the first part of the statement, it suffices to show that when we arrange the vertices of $W$
along $\oC$, for any two consecutive vertices $u$ and $v$ from the arrangement, we have 
$\dist_{\oC}(u,v) \ge t+3$. Note that $V(u^+\oC v^-)\cap W=\emptyset$ for such pairs of $u$ and $v$. 
Assume to the contrary that  there are distinct $u,v\in W$ with $V(u^+\oC v^-)\cap W=\emptyset$  and $\dist_{\oC}(u,v)  \le t+2$. 
Let $C^*=u\iC vv^*Pu^*u$. Since $H$ is complete and $V(u^+\oC v^-)\cap W=\emptyset$,   $H-V(P)$ is a component of $G-V(C^*)$ of order at most $\lambda-1$
and $u^+\oC v^-$ is a component of $G-V(C^*)$
of order at most $t+1$. 
By~\eqref{eqn:degree-on-C}, for each vertex $x\in V(u^+\oC v^-)$, $\deg_G(x, C^*)>\frac{n}{t+1}-1$. 
Applying Lemma~\ref{cycle-extendabilit y}, we find a cycle $C^{**}$ of $G$
with $V(C^{**})=V(C^*)\cup V(u^+\oC v^-)$.  Since $V(G)\setminus V(C^{**})=V(H)\setminus V(P)$,
$C^{**}$ is a $D_\lambda$-cycle of $G$, contradicting the choice of $C$. 

For the second part of the statement, 
we assume to the contrary that $E_G(L^+_u, L^+_v) \ne \emptyset$. 
Applying the first part,   we know that $\dist_{\oC}(u,v) \ge t+3$
and $\dist_{\oC}(v,u) \ge t+3$ (exchanging the role of $u$ and $v$). 
Thus $L^+_u\cap L^+_v=\emptyset$. 
We choose  $x\in L^+_u$ with $\dist_{\oC}(u,x)$ minimum
and $y\in L^+_v$ with $\dist_{\oC}(v,y)$ minimum such that 
$xy\in E(G)$. By this choice of $x$ and $y$,  it follows that 
$E_G(V(u^+ \oC x^-), V(v^+ \oC y^-))=\emptyset$.
Let $C^*=u\iC yx\oC vv^*Pu^*u$. Since $H$ is complete of order $\lambda$ and 
no vertex of $H$ is adjacent in $G$ to any  vertex of $u^+ \oC x^-$ or $v^+ \oC y^-$ by the first part of the statement,   $H-V(P)$ is a component of $G-V(C^*)$ of order at most $\lambda-1$. Also
$u^+ \oC x^-$  and $v^+ \oC y^-$ are components of $G-V(C^*)$ 
of order at most $t+1$.  Since $E_G(V(u^+ \oC x^-), V(v^+ \oC y^-))=\emptyset$, by~\eqref{eqn:degree-on-C}, for each vertex $w\in V(u^+ \oC x^-)\cup V(v^+ \oC y^-)$, $\deg_G(w, C^*)>\frac{n}{t+1}-1$.
Applying Lemma~\ref{cycle-extendabilit y}, we find a cycle $C^{**}$ of $G$
with $V(C^{**})=V(C^*)\cup  V(u^+ \oC x^-)\cup V(v^+ \oC y^-)$.  Since $V(G)\setminus V(C^{**})=V(H)\setminus V(P)$,
$C^{**}$ is a $D_\lambda$-cycle of $G$, contradicting the choice of $C$. 
\qed 

\begin{CLA}\label{claim:W-size}
	$\omega \le \frac{n}{2(t+1)}-\frac{1}{2}$. 
\end{CLA}

\pf Assume otherwise that $\omega > \frac{n}{2(t+1)}-\frac{1}{2}$. 
By Claim~\ref{claim:segments}, for any two distinct $u,v\in W$, 
$G[L_u^+]$ and $G[L_v^+]$ are remote, and $G[L_u^+]$ and $H$
are remote.
Thus  in $G$, there are $\omega+1$ pairwise remote subgraphs. 
By the definition, $G[L_u^+]$ has order $t+2$ for each $u\in W$. 
Let $S=(V(G)\setminus V(H)) \setminus (\bigcup_{u\in W}L_u^+)$. 
Then  $|S| \le  n- \lambda-\frac{n(t+2)}{2(t+1)}+\frac{t+2}{2}  \le \frac{tn}{2(t+1)}+\frac{t}{2}$. Thus 
\begin{eqnarray*}
	\frac{|S|}{c(G-S)} & \le & \frac{\frac{tn}{2(t+1)}+\frac{t}{2}}{\omega+1}<\frac{\frac{tn}{2(t+1)}+\frac{t}{2}}{\frac{n}{2(t+1)}+\frac{1}{2}}=t, 
\end{eqnarray*}
contradicting the toughness of $G$. 
\qed 

Since $\omega \ge 2t$, by Claim~\ref{claim:W-size}, 
we have 
\begin{equation}\label{eqn:n-bound}
	n\ge 4t(t+1). 
\end{equation}
\begin{CLA}\label{claim:H-and-W-size2}
	$\lambda+ \omega  \le \frac{3n}{4(t+1)}+t$. 
\end{CLA}

\pf  Assume to the contrary that $\lambda+ \omega  > \frac{3n}{4(t+1)}+t$.     
By Claim~\ref{claim:W-size}, we know that $\omega \le \frac{n}{2(t+1)}-\frac{1}{2}$.
Since $\omega \ge 2t$,  Lemma~\ref{lem:dominating-cycle} implies that $\lambda \le \frac{n}{2t+1}-t$. 
Thus $\omega >\frac{3n}{4(t+1)}-\frac{n}{2t+1}+2t$. 
By Lemma~\ref{lem:dominating-cycle}, we have 
\begin{eqnarray*}
	n &\ge& (\lambda+t)(\omega+1) \\
	&>&  \left(\frac{3n}{4(t+1)}+t-\omega+t\right) (\omega+1)  \quad \quad \text{($\lambda+ \omega  > \frac{3n}{4(t+1)}+t$ by the assumption.)}\\
	&\ge & \begin{cases}
		\frac{n}{2t+1}(\frac{3n}{4(t+1)}+2t-\frac{n}{2t+1}+1), & \text{if $f(\omega)=\left(\frac{3n}{4(t+1)}+t-\omega+t\right) (\omega+1)$ is increasing;}\\
		(\frac{n}{4(t+1)}+2t+\frac{1}{2})(\frac{n}{2(t+1)}+\frac{1}{2}), & \text{if $f(\omega)=\left(\frac{3n}{4(t+1)}+t-\omega+t\right) (\omega+1)$ is decreasing;}
	\end{cases}\\
	&>& \begin{cases}
		\frac{n}{2t+1}(2t+1)=n,  \quad\quad  \quad\quad\quad\quad\quad\quad\quad\quad \text{if $f(\omega)$ is increasing;}\\
		(\frac{n}{4(t+1)}+2t+\frac{1}{2})\frac{n}{2(t+1)}+\frac{n}{8(t+1)}+t \,\,\, \,\quad  \quad\text{if $f(\omega)$ is decreasing;}\\
		>
		\begin{cases}
			(2t+2t+\frac{1}{2})\frac{n}{2(t+1)}>n,  &\text{if $n\ge 8t(t+1)$}; \\ 
			(3t+\frac{1}{2})\frac{n}{2(t+1)}+\frac{n}{8(t+1)}+t=\frac{3tn}{2(t+1)}+\frac{3n+8t(t+1)}{8(t+1)}>n, &\text{if $n< 8t(t+1)$};
		\end{cases}
	\end{cases}
\end{eqnarray*}
achieving a contradiction, where $n\ge 4t(t+1)$ was used to obtain $\frac{n}{4(t+1)}\ge t$  in the last inequality when $f(\omega)$ is decreasing.  
\qed 

By  Claim~\ref{claim:one-component} and Claim~\ref{claim:H-and-W-size2}, we have 
\begin{equation}\label{eqn:degree2}
	\deg_G(v,C) >\frac{1.25n}{t+1}-1 \quad\text{for any $v\in V(C)\setminus W$}. 
\end{equation}

We will now explore  the neighborhood of vertices from $W^+:=\{w\in V(C): w^-\in W\}$,
and show that some vertices from the neighborhood have similar properties as those 
in $W^+$. 
By Claim~\ref{claim:t>=1}, we know that $|W| \ge 3$ and so $|W^+|\ge 3$. 
Equation~\eqref{eqn:degree2} allows us to construct the vertex sets $L_x^-(\frac{0.25n}{t+1}+2)$ for each $x\in N_C(W^+)$. 
For notation simplicity,  we use $L_x^-$ for $L_x^-(\frac{0.25n}{t+1}+2)$.  Note that the statement below is not true in general if 
we replace $L_x^-(\frac{0.25n}{t+1}+2)$ by $L_x^+(\frac{0.25n}{t+1}+2)$. 
\begin{CLA}\label{claim:crossing-chords}
	Let $u\in W^+$ and $x\in N_C(u)$.   Then 
	\begin{enumerate}[(1)]
		\item $L_x^-\cap W=\emptyset$. 
		\item Let $v\in W^+$ and $y\in N_C(v)$ such that $ux$ and $ vy$ are two crossing chords of $C$. Then $\dist_{\oC}(x,y) \ge \frac{0.25n}{t+1}+3$. 
	\end{enumerate}
\end{CLA}

\pf 
For Statement (1), 
suppose to the contrary that there exists  $z\in W$ such that  $z\in L_x^-$. Then $\dist_{\oC}(z,x)  \le \frac{0.25n}{t+1}+2$. 
We  choose $z\in W$ with $\dist_{\oC}(z,x) $ minimum.  Then $V(z^+\oC x^-)\cap W=\emptyset$ and $\dist_{\oC}(z,x)  \le \frac{0.25n}{t+1}+2$. 
Let $u^* \in N_H(u^-)$, $z^* \in N_H(z)$, and $P^*$ be a $(u^*,z^*)$-path of $H$.
Then $C^*=z\iC ux \oC u^-u^*P^*z^*z$ is a cycle. 
As $H$ is complete of order $\lambda$ and $ V(z^+\oC x^-)\cap W=\emptyset$, 
we know that $H-V(P^*)$
is a component of $G-V(C^*)$ of order at most $\lambda-1$. Also, $z^+\oC x^-$ is a component of 
$G-V(C^*)$ of order at most $\frac{0.25n}{t+1}+1$.  By~\eqref{eqn:degree2}, for each vertex $w\in V(z^+\oC x^-)$, $\deg_G(w, C^*)>\frac{n}{t+1}-1$. 
Applying Lemma~\ref{cycle-extendabilit y}, we find a cycle $C^{**}$ of $G$
with $V(C^{**})=V(C^*)\cup V(z^+\oC x^-)$.  Since $V(G)\setminus V(C^{**})=V(H)\setminus V(P^*)$,
$C^{**}$ is a $D_\lambda$-cycle of $G$, contradicting the choice of $C$.  

Let $u^*\in N_H(u^-), v^*\in N_H(v^-)$ and $P$ be a $(u^*,v^*)$-path of $H$. 
For Statement (2), suppose to the contrary that  $\dist_{\oC}(x,y)  \le \frac{0.25n}{t+1}+2$.  
We assume without loss of generality that $u,v,x,y$ appear in this order along $\oC$. 
Let   $C^*=u\oC v^-v^*Pu^*u^-\iC yv \oC xu$.
Since $H$ is complete of order $\lambda$ and $V(x^+\oC y^-)\cap W=\emptyset$ by Statement (1) (note $V(x^+\oC y^-) \subseteq L_y^-$),   $H-V(P)$ is a component of $G-V(C^*)$
of order at most $\lambda-1$.
Also, $x^+\oC y^-$ is a  component of $G-V(C^*)$
of order at most $\frac{0.25n}{t+1}+1$. 
By~\eqref{eqn:degree2}, for each vertex $w\in V(x^+\oC y^-)$, $\deg_G(w, C^*)>\frac{n}{t+1}-1$. 
Applying Lemma~\ref{cycle-extendabilit y}, we find a cycle $C^{**}$ of $G$
with $V(C^{**})=V(C^*)\cup V(x^+\oC y^-)$.  Since $V(G)\setminus V(C^{**})=V(H)\setminus V(P)$,
$C^{**}$ is a $D_\lambda$-cycle of $G$, contradicting the choice of $C$. 
\qed

For two distinct vertices $x,y\in N_C(W^+)$, we say $x$
and $y$ form a \emph{crossing}  if there exist distinct vertices $u, v\in W^+$
such that $ux$ and $vy$ are crossing chords of $C$.   
By Claim~\ref{claim:crossing-chords}(2), there are at least $\frac{0.5n}{t+1}+2$ vertices between $x$
and $y$  along $\oC$ for any two $x,y\in N_C(W^+)$ such that $x$ and $y$ form a crossing. 
Our goal below is to find at least $\frac{n}{2(t+1)}$ vertices from $N_C(W^+)$ such that there are at least  $\frac{0.5n}{t+1}+2$ vertices 
between any two of them along $\oC$. Then we will reach a contradiction by showing that $|V(C)|\ge n$.  
Define 
\begin{equation}\nonumber
	A=\{u\in V(C)\,:\, \deg_G(u, W^+) =1\} \quad  \text{and} \quad B=\{u\in V(C)\,:\, \deg_G(u, W^+) \ge 2\}.
\end{equation}
Let $u\in W^+$ and $p=\deg_G(u, B)$ for some positive integer $p$, and let $$N_G(u)\cap B=\{x_1,x_2, x_3, \ldots, x_p\}.$$  We may assume that $x_1, x_2,\ldots, x_p$  appear in  the same order  along $\oC$.  We separate those vertices according to vertices of $W$. 
By Claim~\ref{claim:crossing-chords}(1), we have $L_{x_i}^-\cap W=\emptyset$ for each $i\in [1,p]$. 
Therefore for some integer $q\ge 1$, we assume that  $x_1,\ldots ,x_p$ are grouped into $q$ sets
$$
B_1=\{x_{b_0+1}, \ldots, x_{b_1}\}, \quad  B_2=\{x_{b_1+1}, \ldots, x_{b_2}\}, \quad \ldots, \quad B_q=\{x_{b_{q-1}+1}, \ldots, x_{b_q}\}, 
$$
where $b_0=0$ and $b_q=p$, 
such that  $V(x_{b_j+1}^+\oC x_{b_{j+1}}^-)\cap W=\emptyset$ for each $j\in[0,q-1]$.   Furthermore, we may assume that the number $q$ of sets
with the property above is minimum. As $W \ne \emptyset$, the minimality of $q$ in turn implies $V(x_{b_j} \oC x_{b_j+1}) \cap W \ne \emptyset$ for each 
$j\in [1,q]$, where $x_{b_q+1}:=x_1$. Hence, by Claim~\ref{claim:crossing-chords}(1), we have 
\begin{equation}\label{eqn:W}
	\dist_{\oC}(x_{b_j}, x_{b_j+1}) \ge \frac{0.5n}{t+1}+3. 
\end{equation}

\begin{CLA}\label{claim:half}
	For each $i\in [1,q]$, $B_i$ has at least $|B_i|/2$ vertices such that the distance between any two of them on $C$ is at least $\frac{0.5n}{t+1}+3$. 
\end{CLA}

\pf 
We partition  $B_i$
into two subsets according to whether or not vertices in $N_G(x)\cap (W^+\setminus \{u\})$ fall into $x_{b_i}\oC u$ for $x\in B_i$. 
Define 
\begin{eqnarray*}
	B_{i1}=\{x \in B_i\,:\,   (N_G(x)\cap (W^+\setminus \{u\}))\cap V(x_{b_i}\oC u) \ne \emptyset\}, &  B_{i2}=B_i\setminus B_{i1}.
\end{eqnarray*}
By the Pigeonhole Principle, we have $|B_{i1}| \ge |B_i|/2$ or $|B_{i2}| \ge |B_i|/2$. 
We show that any two distinct vertices from $B_{i1}$ or from $B_{i2}$ have distance at least $\frac{0.5n}{t+1}+3$ between them on $C$. 
Let $x_a,x_b\in B_{i1}$ or $x_a,x_b\in B_{i2}$ be distinct, where $a,b\in [b_{i-1}+1,b_i]$. 
If $a>b$, then $\dist_{\oC}(x_a, x_b) \ge \dist_{\oC}(x_{b_i}, x_{b_{i-1}+1}) \ge \dist_{\oC}(x_{b_i}, x_{b_i+1})$. As $\dist_{\oC}(x_{b_i}, x_{b_i+1}) \ge \frac{0.5n}{t+1}+3$ by
~\eqref{eqn:W}, we have $\dist_{\oC}(x_a, x_b) \ge \frac{0.5n}{t+1}+3$. 
Thus we assume $a<b$. 

If $x_a,x_b\in B_{i1}$,
since $V(x_{b_{i-1}+1}^+\oC x_{b_{i+1}}^-)\cap W=\emptyset$,  then we know $x_{b_i} \not\in W^+$. 
Thus by the definition of $B_{i1}$, there exists $v\in (N_G(x_a)\cap (W^+\setminus \{u\}))\cap V(x_{b_i}^+\oC u)$. Then the four vertices $u,v, x_a,  x_b$ appear in the order $x_a, x_b, v, u$ along $\oC$ and so $vx_a$ and $ux_b$ are crossing chords of $C$. By Claim~\ref{claim:crossing-chords}(2), we have $\dist_{\oC}(x_a,x_b) \ge \frac{0.5n}{t+1}+3$. 

Suppose now that $x_a,x_b\in B_{i2}$. If $a>b_{i-1}+1$ or $a=b_{i-1}+1$ but $x_a \not\in W$, then $W^+\cap V(x_a\oC x_b)=\emptyset$
by the property of $B_i$  that   $V(x_{b_{i-1}+1}^+\oC x_{b_{i+1}}^-)\cap W=\emptyset$.  By the definition of $B_{i2}$, there exists $v\in (N_G(x_b)\cap (W^+\setminus \{u\}))\cap V( u\oC x_{b_{i-1}+1})$. Then the four vertices $u, v, x_a, x_b$ appear in the order $ x_a, x_b, u,v$ along $\oC$ and so $ux_a$ and $vx_b$ are  crossing chords of $C$. By Claim~\ref{claim:crossing-chords}(2), we have $\dist_{\oC}(x_a,x_b) \ge \frac{0.5n}{t+1}+3$. 
Thus we assume  $a=b_{i-1}+1$ and $x_a \in W$. By Claim~\ref{claim:crossing-chords}(1), we know that $\dist_{\oC}(x_a,x_{a+1}) \ge \frac{0.5n}{t+1}+3$
and so $\dist_{\oC}(x_a,x_b) \ge \frac{0.5n}{t+1}+3$. 
\qed

By Claim~\ref{claim:half}, for each $i\in [1,q]$, we take a subset of at least $|B_i|/2$ vertices from $B_i$ such that the distance between any two of 
them on $C$ is at least $\frac{0.5n}{t+1}+3$. We let $\{ y_1, y_2, \ldots, y_k\}$ be the union of all these $q$ subsets of vertices. Then $k\ge \lceil \frac{p}{2} \rceil$. We further assume 
those vertices appear in the order $y_1, \ldots, y_k$ along $\oC$. 
For any two distinct $y_j, y_\ell$ with $j,\ell \in [1,k]$, if $y_j, y_\ell$ are from the same $B_i$ for some $i\in [1,q]$, then we
have $\dist_{\oC}(y_j,y_\ell) \ge \frac{0.5n}{t+1}+3$ by Claim~\ref{claim:half}. Otherwise, by~\eqref{eqn:W}, we also have $\dist_{\oC}(y_j,y_\ell) \ge \frac{0.5n}{t+1}+3$. 
Thus by~\eqref{eqn:n-bound} that $n\ge 4t(t+1)$, we have 
$
n>|V(C)| \ge k \left(\frac{0.5n}{t+1}+2 \right) \ge k(2t+2). 
$
This inequality  implies  $k\le \frac{n}{2(t+1)}$. Therefore 
\begin{equation}\label{eqn:u-degree}
	\deg_G(u,B) \le 2k \le \frac{n}{t+1} \quad \text{for any $u\in W^+$}. 
\end{equation}

Let $s=\sum\limits_{v\in W^+} \deg_G(v, B)$ for some positive integer $s$. Then there exists $u\in W^+$
such that $\deg_G(u,B) \ge \frac{s}{|W^+|}=\frac{s}{\omega}$. Following the notation defined above, we let 
$\{y_1, \ldots, y_k\} \subseteq N_G(u)\cap B$ such that $\dist_{\oC}(y_i,y_j) \ge \frac{0.5n}{t+1}+3$ 
for any  distinct $i, j\in [1,k]$, where  $k\ge \frac{1}{2}\deg_G(u,B)$. 
By Claim~\ref{claim:crossing-chords}(1), we have $L_{y_i}^-\cap W=\emptyset$
for any $i\in [1,k]$. 
Thus, as $n\ge 4t(t+1)$ by~\eqref{eqn:n-bound} and $t>1$, we have $\dist_{\oC}(y_i,w^+) \ge \frac{0.5n}{t+1}+2 \ge 2t+2>t+3$ for  any $w\in W$. 
Thus for any $i \in [1,k]$,  $y_i^+ \oC y_{i+1}^-$ has  at least $t+2$
vertices that are nonadjacent in $G$ to any vertex of $W^+$, where $y_{k+1}:=y_1$.  
Hence 
\begin{eqnarray*}
	|A|&\le& n-|B|-\frac{1}{2}\deg_G(u,B)(t+2) \\
	& \le& n-\frac{s}{\omega}-\frac{s}{2\omega}(t+2),
\end{eqnarray*}
since $|B| \ge \frac{s}{\omega}$.   
By~\eqref{eqn:degree2}, we get 
\begin{equation}\label{eqn:sumdegree}
	\omega\left( \frac{1.25n}{t+1}-1\right)<\sum\limits_{u\in W^+} d_G(u) =|A|+s \le n-\frac{s}{\omega}-\frac{s}{2 \omega}(t+2)+s.
\end{equation}
Since $\omega \ge 2t$, $n\ge 4t(t+1)$ by~\eqref{eqn:n-bound} and $t>1$,
it follows that $\omega( \frac{1.25n}{t+1}-1) \ge \frac{2tn}{t+1}>n$.
Thus, $-\frac{s}{\omega}-\frac{s}{2 \omega}(t+2)+s>0$
and so $2\omega-t-4>0$. 
Thus by~\eqref{eqn:sumdegree}, 
$$
s> \frac{2\omega^2(\frac{1.25n}{t+1}-1)-2n\omega}{2\omega-t-4}.
$$	 
Next, we claim 
\begin{equation}\label{eqn:sumdegree2}
	\frac{2\omega^2(\frac{1.25n}{t+1}-1)-2n\omega}{2\omega-t-4}>\frac{n\omega}{t+1}, 
\end{equation}
which will in turn give $s>\frac{n\omega}{t+1}$. 
As $s=\sum\limits_{v\in W^+} d_G(v, B)$ and $|W^+|=\omega$, it then will follow that there exists $u\in W^+$ with $\deg_G(u, B)>\frac{n}{t+1}$, and so will give a contradiction to~\eqref{eqn:u-degree}. 
To prove~\eqref{eqn:sumdegree2}, it suffices to show that 
$2\omega(\frac{1.25n}{t+1}-1)-2n>\frac{n(2\omega-t-4)}{t+1}$, which is true as shown below:
\begin{eqnarray*}
	\frac{2n}{t+1}-4t&>&0 \,\, \text{ (since $n \ge 4t(t+1)$ by~\eqref{eqn:n-bound})} \quad \Rightarrow \\
	\frac{0.5\omega n}{t+1}+\frac{(t+4)n}{t+1}-2\omega-2n&>&0\,\,\text{(the function on the left  increases in $\omega$, $\omega \ge 2t$)} \quad \Leftrightarrow\\
	2\omega(\frac{1.25n}{t+1}-1)-2n&>&\frac{n(2\omega-t-4)}{t+1}. 
\end{eqnarray*}
\qed

\section*{Acknowledgments}

The author is very grateful to the two anonymous referees  for their careful reading and valuable
comments, which greatly  improved this paper.

\end{document}